\documentclass{amsart}
\usepackage{amsmath}
\usepackage{amsthm}
\usepackage{amssymb}
\usepackage{epsfig}

\theoremstyle{plain}
\newtheorem{thm}{Theorem}[section]
\newtheorem{cor}[thm]{Corollary}

\newtheorem{prop}[thm]{Proposition}
\theoremstyle{definition}
\newtheorem{defn}[thm]{Definition}
\theoremstyle{remark}

\title{Spotlight Tiling}
\author{Bridget Eileen Tenner}
\address{Department of Mathematical Sciences, DePaul University, 2320 North Kenmore Avenue, Chicago, IL 60614, USA}
\email{bridget@math.depaul.edu}
\date{June 25, 2008}

\begin{document}

\begin{abstract}
This article introduces spotlight tiling, a type of covering which is similar to tiling.  The distinguishing aspects of spotlight tiling are that the ``tiles'' have elastic size, and that the order of placement is significant.  Spotlight tilings are decompositions, or coverings, and can be considered dynamic as compared to typical static tiling methods.  A thorough examination of spotlight tilings of rectangles is presented, including the distribution of such tilings according to size, and how the directions of the spotlights themselves are distributed.  The spotlight tilings of several other regions are studied, and suggest that further analysis of spotlight tilings will continue to yield elegant results and enumerations.
\end{abstract}

\maketitle

\section{Introduction}

Domino tilings, and relatedly perfect matchings, are well studied objects in combinatorics and statistical mechanics.  In the typical setup, there is a finite set $S$ of distinct tiles which may be used repeatedly to tile a particular region or family of regions.  It is then natural to count the number of ways a particular region can be tiled by elements of $S$, or, more fundamentally, to determine if any such tiling is even possible.  The number of domino tilings of a rectangle, the most elementary region, was computed by Kasteleyn in \cite{kasteleyn}.  

The number of tilings of an $m \times n$ rectangle can become much simpler if certain restrictions are imposed.  For example, suppose that the region $R$ is colored as a checkerboard having a black upper-left square, with alternating black and white squares in each column or row.  Restrict the set $S$ to contain vertical dominos of both colorings (one with a white top square and one with a black top square), and only the horizontal domino with a black left square.  Then, it is straightforward to show that the number of such tilings of an $m \times n$ region $R$ by elements of $S$ is

\begin{equation*}
\left\{\begin{array}{c@{\quad:\quad}l}
0\phantom{^{n/2}}  & m \text{ and } n \text{ are both odd};\\
1\phantom{^{n/2}}  & m \text{ is even};\\
\left(\frac{m+1}{2}\right)^{n/2} & m \text{ is odd and } n \text{ is even}.
\end{array}\right.
\end{equation*}

\noindent These numbers are sequence A133300 of \cite{oeis}.  There is a rich literature concerning domino tilings, as well as tilings by shapes which are generalizations of dominoes in some aspect.  For example, see \cite{golomb, kasteleyn, kenyon, propp}.

Typical tiling results do not depend on the order in which the tiles are placed.  Because the set $S$ of allowable tiles does not change as each tile is placed in the region, tiles may be considered to be placed simultaneously.

This article introduces a method of covering regions, somewhat related to tilings, and provides a sample of results answering the most basic questions about this method.  There are two significant differences between this and previous tiling methods: the shape of the ``tiles'' here is elastic, and the order in which they are positioned is important.  One interpretation of these differences is that the method studied here is a \emph{dynamic} covering model, while other methods, such as domino tiling, would be static.

Henceforth, the ``tiles'' in this paper will be called \emph{spotlights} to emphasize their elastic nature and to avoid confusion with more customary notions of tiling.
  
In this initial foray into the dynamic spotlight tiling model, the rules for placing the spotlights will be somewhat strict, requiring that each spotlight originate in the same type of corner.  Relaxing this restriction leads to other interesting questions, discussed in the last section of the paper.

As mentioned above, spotlights are placed in the region sequentially, and after each placement the set of allowable spotlights may change.  To be specific, first a particular corner direction is specified (\emph{northwest} for the duration of this article).  At each stage a spotlight is placed with one end point in a ``corner,'' as defined by the chosen direction, and the spotlight must extend as far as possible from this corner either horizontally or vertically.  This type of covering is called a \emph{spotlight tiling}, in reference to the fact that it is like placing a spotlight in one of the specified corners and turning it to point horizontally or vertically so that it shines as far as possible until it reaches an obstruction.

Spotlight tilings of rectangles are examined thoroughly below, including a description of various statistics, such as the number of spotlights needed and the average number of spotlights used in a spotlight tiling of the rectangle.  Additionally, spotlight tilings of regions which are similar to rectangles are studied.  The nature of spotlight tiling means that many of the proofs used to obtain the results below are recursive in nature.

The most basic region is an $m \times n$ rectangle.  Therefore, in this introductory analysis of spotlight tiling, attention is primarily focused on rectangles, in terms of their enumeration and their properties.  This will be the substance of Section~\ref{sec:tiling rectangles}.  For example, in addition to determining the number of spotlight tilings of an $m \times n$ rectangle, more detailed statistics will be studied.  Unlike other sorts of tilings, where the number of tiles required to cover a region is fixed, the number of spotlights used depends on the particular spotlight tiling itself.  The distribution of the number of these tiles will be part of the discussion in Section~\ref{sec:tiling rectangles}.  Following this discussion, in Section~\ref{sec:other regions}, attention will be turned to spotlight tilings of regions which are formed from rectangles by removing squares at the corners.  The recursive nature of these spotlight tilings leads naturally to recursive enumeration formulae.  In some cases, these equations will be left in a recursive format, as it is simpler to read them in this manner.  In other situations, when a closed form itself is quite elegant, both the recursive and the closed formulae will be given.  Finally, in Section~\ref{sec:frames}, the spotlight tilings of a certain family of frame-like regions is explored.  The paper concludes with a brief discussion of how spotlight tilings may be studied further.

\section{Definitions}

The basic definitions and notation of this article are outlined below.

\begin{defn}
A \emph{region} is the dual of a finite connected induced subgraph of $\mathbb{Z}^2$.
\end{defn}

Spotlight tilings rely on the choice of a particular direction and type of corner, in this case a northwest corner.

\begin{defn}
A \emph{northwest corner} in a region is a square belonging to the region that is bound above and on the left by the boundary edge of the region.
\end{defn}

\begin{figure}[htbp]
\epsfig{file=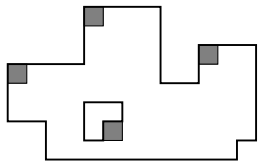,scale=.4}
\caption{A region with four northwest corners, which are marked by shading.}\label{fig:nwcorners}
\end{figure}

For example, the four northwest corners of the region in Figure~\ref{fig:nwcorners} have been shaded.

As discussed in the introduction, spotlight tilings differ in nature from static tilings.  Instead of choosing from a finite set of tiles, the possible spotlights themselves are defined by the region and any spotlights that have been positioned previously.

\begin{defn}\label{defn:spotlight tile}
A \emph{spotlight} with an endpoint in a northwest corner $s$ extends as far east horizontally or south vertically from $s$ as possible, terminating at the boundary of the region, or when it encounters a spotlight that has already been placed.
\end{defn}

\begin{defn}\label{defn:spotlight tiling}
Given a region $R$, a \emph{spotlight tiling} of $R$ is defined recursively as follows.  Choose any northwest corner $s \in R$.  Place a spotlight tile with an endpoint in $s$, extending either horizontally (east) or vertically (south) as far as possible.  Let $R'$ be the collection of disjoint regions remaining after placing this spotlight in $R$.  The spotlight tiling of $R$ is completed by finding spotlight tilings of each connected component of $R'$.
\end{defn}

\begin{figure}[htbp]
\parbox{.8in}{\epsfig{file=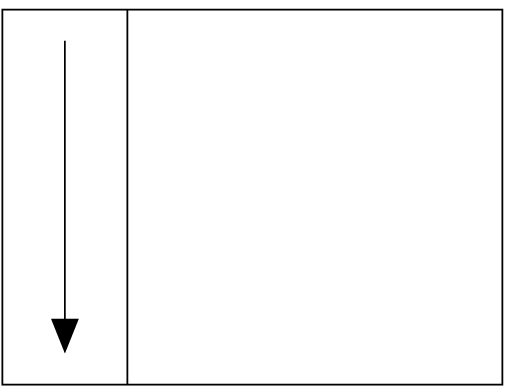,scale=.4}} $\Rightarrow$ \parbox{.8in}{\epsfig{file=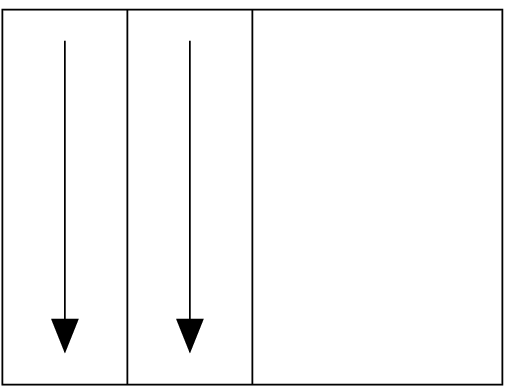,scale=.4}} $\Rightarrow$ \parbox{.8in}{\epsfig{file=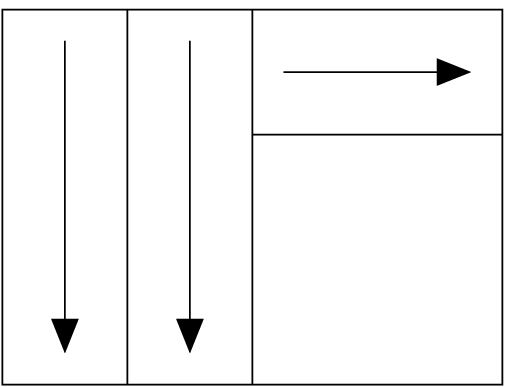,scale=.4}} $\Rightarrow$ \parbox{.8in}{\epsfig{file=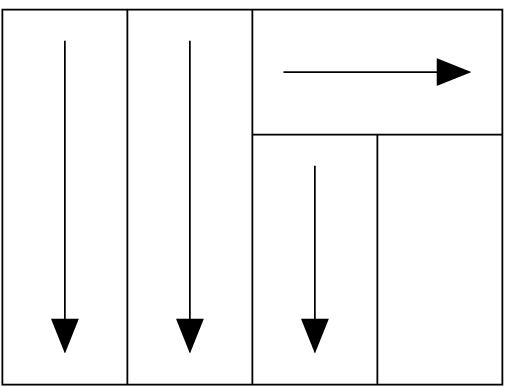,scale=.4}} $\Rightarrow$ \parbox{.8in}{\epsfig{file=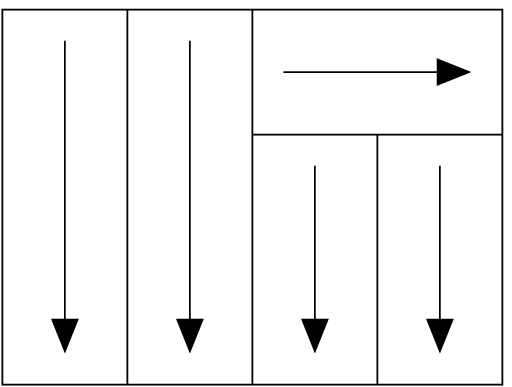,scale=.4}}
\caption{The recursive construction of a spotlight tiling of a $3 \times 4$ rectangle.  The arrows are provided here only to highlight the direction (horizontal or vertical) of each spotlight.} \label{fig:3x4}
\end{figure}

A spotlight tiling of a $3 \times 4$ rectangle is depicted in Figure~\ref{fig:3x4}.  The complete tiling is the last image in the figure, having been built successfully from the previous images.

Although spotlight tiles are placed sequentially in a region, two spotlight tilings are considered distinct only if they look different once all the spotlights are in place.  In other words, if there is more than one order in which the spotlights can be placed in the region, this alone does not distinguish one tiling from another.  Moreover, the direction (horizontal or vertical) of a spotlight is obvious except in certain cases of tiles of length one, where the direction of such a spotlight will not be specified as uniquely horizontal or vertical.  Ignorance of the orientation of this spotlight maintains consistency with the fact that two spotlight tilings differ only if they look different.  However, the enumerations of this paper could be reformulated without this stipulation, and similarly nice results would ensue.

The order in which spotlights are placed in a spotlight tiling of a region $R$ can be recovered in some cases.  More precisely, a complete recovery is possible if the region $R$ has only one northwest corner and does not have any holes.  If $R$ did have holes, then it could be possible to place some number of spotlights in $R$ and yield an untiled subregion having more than one northwest corner.

\begin{defn}
If the last spotlight placed in a spotlight tiling has length $1$, it is a \emph{HV-spotlight}, referring to the fact that the spotlight's direction could be considered to be either horizontal or vertical.
\end{defn}

\begin{figure}[htbp]
\epsfig{file=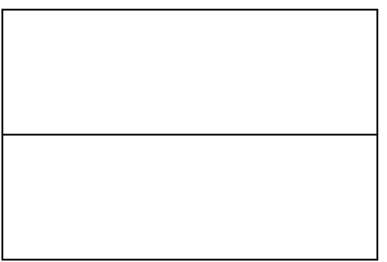,scale=.4} \hspace{.25in} \epsfig{file=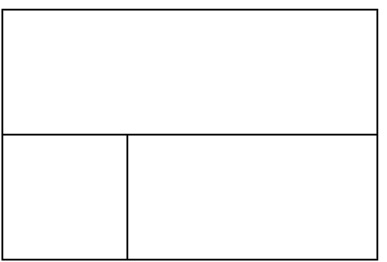,scale=.4} \hspace{.25in} \epsfig{file=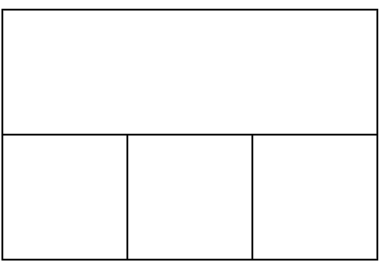,scale=.4} \hspace{.25in} \epsfig{file=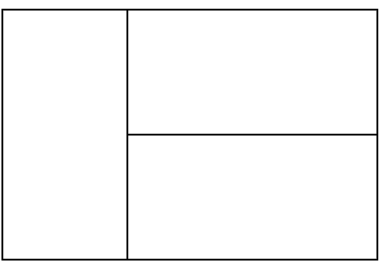,scale=.4}\\
\vspace{.25in}
\epsfig{file=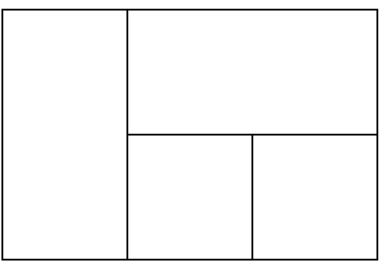,scale=.4} \hspace{.25in} \epsfig{file=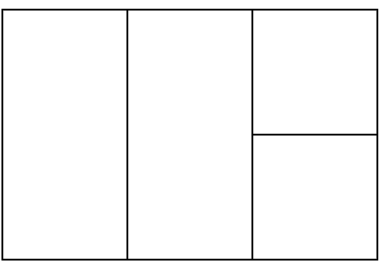,scale=.4} \hspace{.25in} \epsfig{file=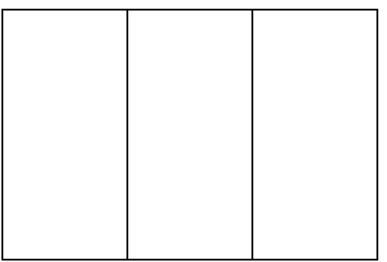,scale=.4}
\caption{The seven distinct spotlight tilings of a $2 \times 3$ rectangle.  In the third, fifth, and sixth of these, the last (southeast-most) spotlight is a HV-spotlight.}\label{fig:2x3}
\end{figure}

The seven different spotlight tilings of a $2 \times 3$ rectangle are depicted in Figure~\ref{fig:2x3}.

\begin{defn}
Let $R_{m,n}$ denote an $m \times n$ rectangle.  The set of spotlight tilings of $R_{m,n}$ is denoted $\mathcal{T}_{m,n}$, and $T_{m,n} = |\mathcal{T}_{m,n}|$.  For all $m,n > 0$, set $T_{m,0} = T_{0,n} = 1$.  

\end{defn}

As depicted in Figure~\ref{fig:2x3}, $T_{2,3} = 7$.

The recursive definition of spotlight tiling means that
\begin{equation}\label{eqn:set recursion}
\begin{split}
\mathcal{T}_{m,n} =& \left\{\text{one } (1 \times n)\text{-spotlight together with } t \mid t \in \mathcal{T}_{m-1,n}\right\}\\
&\cup \left\{\text{one } (m \times 1)\text{-spotlight together with } t \mid t \in \mathcal{T}_{m,n-1}\right\}.
\end{split}
\end{equation}

\section{Spotlight tilings of rectangles}\label{sec:tiling rectangles}

The first goal of this examination of spotlight tilings is a thorough understanding of spotlight tilings of rectangles.  Since the definition of a spotlight tiling gives no preference to horizontal or vertical spotlights, all results in this section should be symmetric with respect to $m$ and $n$.  In particular, it should be the case that $T_{m,n} = T_{n,m}$.

A precise formula for $T_{m,n}$ is straightforward to compute, based on the recursive nature of Definition~\ref{defn:spotlight tiling}.

\begin{thm}\label{thm:total rect}
For all $m,n \ge 1$,
\begin{equation}\label{eqn:rectangle value}
T_{m,n} = \binom{m+n}{m} - \binom{m+n-2}{m-1}.
\end{equation}
\end{thm}

\begin{proof}
Definition~\ref{defn:spotlight tiling} gives the recursive formula
\begin{equation}\label{eqn:rectangle recursion}
T_{m,n} = T_{m-1,n} + T_{m,n-1}
\end{equation}
\noindent for all positive $m$ and $n$ such that $mn >1$.  Since $T_{1,1} = 1$, equation~\eqref{eqn:rectangle value} is satisfied for $m=n=1$.  Supposing inductively that the result holds whenever the dimensions of the rectangle sum to less than $k$, consider an $m \times n$ rectangle where $m + n = k$. Then, using equation~\eqref{eqn:rectangle recursion},
\begin{eqnarray*}
T_{m,n} &=& T_{m-1,n} + T_{m,n-1}\\
&=& \binom{m + n - 1}{m-1} - \binom{m + n - 3}{m - 2} + \binom{m + n - 1}{m} - \binom{m + n - 3}{m - 1}\\
&=& \binom{m + n}{m} - \binom{m + n - 2}{m - 1},
\end{eqnarray*}
\noindent Thus the result holds for all $m, n \ge 1$.
\end{proof}

Notice that equation~\eqref{eqn:rectangle value} is symmetric in $m$ and $n$, as required.  The values of $T_{m,n}$ for small $m$ and $n$ are displayed in Table~\ref{table:T_{m,n}}.  Additionally, these are sequence A051597 of \cite{oeis}.

\begin{table}[htbp]
\centering
\begin{tabular}{c|ccccccc}
\rule[-2mm]{0mm}{6mm}$T_{m,n}$ & $n=1$ & 2 & 3 & 4 & 5 & 6 & 7\\
\hline
\rule[0mm]{0mm}{4mm}$m=1$ & 1 & 2 & 3 & 4 & 5 & 6 & 7\\
2 & 2 & 4 & 7 & 11 & 16 & 22 & 29\\
3 & 3 & 7 & 14 & 25 & 41 & 63 & 92\\
4 & 4 & 11 & 25 & 50 & 91 & 154 & 246\\
5 & 5 & 16 & 41 & 91 & 182 & 336 & 582\\
6 & 6 & 22 & 63 & 154 & 336 & 672 & 1254\\
7 & 7 & 29 & 92 & 246 & 582 & 1254 & 2508\\
\end{tabular}
\smallskip
\caption{The number of spotlight tilings of $R_{m,n}$, for $m,n \in [1,7]$.}\label{table:T_{m,n}}
\end{table}

As demonstrated in Figure~\ref{fig:2x3}, the number of spotlights in a particular spotlight tiling of $R_{m,n}$ is not fixed.  For example, a spotlight tiling of $R_{2,3}$ can consist of $2$, $3$, or $4$ spotlights.  Therefore, to better understand spotlight tilings of rectangles, it is important to understand how many spotlights may (likewise, ``must'' and ``can'') be used in a spotlight tiling of $R_{m,n}$, and how many spotlight tilings of the rectangle use exactly $r$ spotlights.  There are additional aspects of spotlight tilings using the minimal or maximal number of spotlights that are of interest as well.

\begin{defn}
For a spotlight tiling $t$ of a region $R$, let $|t|$ be the number of spotlights used in $t$, known as the \emph{size} of $t$.
\end{defn}

\begin{defn}
Let $t^-_{m,n}$ denote the minimum number of spotlights needed in a spotlight tiling of $R_{m,n}$, and let $t^+_{m,n}$ denote the maximum number of spotlights that can be used in a spotlight tiling of $R_{m,n}$.  That is,
\begin{eqnarray*}
t^-_{m,n} = \min_{t \in \mathcal{T}_{m,n}} |t|\\
t^+_{m,n} = \max_{t \in \mathcal{T}_{m,n}} |t|\\
\end{eqnarray*}
An element of $\mathcal{T}_{m,n}$ using $t^-_{m,n}$ spotlights is a \emph{minimal} spotlight tiling, while one that uses $t^+_{m,n}$ spotlights is a \emph{maximal} spotlight tiling.
\end{defn}

\begin{prop}\label{prop:bounds}
For all $m, n \ge 1$,
\begin{eqnarray}
t^-_{m,n} &=& \min\{m,n\};\label{eqn:t^-}\\
t^+_{m,n} &=& m + n - 1.\label{eqn:t^+}
\end{eqnarray}
\end{prop}

\begin{proof}

By the definition of spotlight tilings, it is clear that the minimum number of spotlights needed depends on the minimum dimension of $R_{m,n}$.  Suppose, without loss of generality, that $m \le n$.  If fewer than $m$ spotlights are placed in $R_{m,n}$, then at least one row and at least one column are not completely covered.  Thus, $t^-_{m,n}$ can be no less than $m$.  Additionally, one spotlight tiling of the rectangle consists of $m$ horizontal spotlights, so $t^-_{m,n} = m$.  This proves equation~\eqref{eqn:t^-}.

Equation~\eqref{eqn:set recursion} implies that $t^+_{m,n} = \max\{1 + t^+_{m-1,n}, 1 + t^+_{m,n-1}\}$.  Then, since $t^+_{1,1} = 1$ and $t^+_{m,1} = m$, the rest of the proof of equation~\eqref{eqn:t^+} follows inductively.
\end{proof}

Note that $t^-_{m,n} = t^+_{m,n}$ if and only if $m = n = 1$.  Therefore, in anything larger than a $1\times 1$ square, there will be variation in the number of spotlights used.

The number of minimal spotlight tilings of an $m \times n$ rectangle is necessarily $1$ or $2$, depending on whether $m \neq n$ or $m = n$.  This will be included in a more general argument in Theorem~\ref{thm:counting rect}.

On the other hand, the number of maximal spotlight tilings is somewhat specialized and will first be treated independently.

\begin{thm}\label{thm:max rect}
The number of maximal spotlight tilings of $R_{m,n}$ is
\begin{equation*}
\binom{m + n - 2}{m - 1}.
\end{equation*}
\end{thm}

\begin{proof}

Equations~\eqref{eqn:set recursion} and~\eqref{eqn:t^+} imply that once the first spotlight has been placed in the rectangle, this can (and, in fact, must) be completed to a maximal tiling of the rectangle by finding a maximal spotlight tiling of the resulting sub-rectangle ($R_{m-1,n}$ or $R_{m,n-1}$, depending on whether the first spotlight was horizontal or vertical).

There is a single element in the set $\mathcal{T}_{1,1}$, and it consists of a single HV-spotlight.  Therefore, using equation~\eqref{eqn:set recursion}, the last spotlight placed in a maximal spotlight tiling must be an HV-spotlight.  In fact, if $m$ and $n$ are not both equal to $1$, then the penultimate spotlight placed in a maximal spotlight tiling of $R_{m,n}$ must also have length $1$, although this will not be an HV-spotlight since its direction must be specified.

The result follows immediately by induction.
\end{proof}

Alternatively, Theorem~\ref{thm:max rect} can also be proved bijectively in the following manner.  By nature of spotlight tiling, there cannot be more than $m$ horizontal spotlights or $n$ vertical spotlights in an element of $\mathcal{T}_{m,n}$.  If the last spotlight is an HV-spotlight, than of the previous $m + n - 2$ spotlights in a maximal spotlight tiling, at most $m - 1$ can be horizontal and at most $n - 1$ can be vertical.  Consequently, of these $m + n - 2$ spotlights, exactly $m - 1$ are horizontal and exactly $n - 1$ are vertical.  Consider an initial set of spotlights in $R_{m,n}$, consisting of at most $m-1$ horizontal spotlights and at most $n-1$ vertical spotlights.  Any such initial spotlight tiling can be completed to a maximal spotlight tiling.  Therefore the number of maximal spotlight tilings depends only on which $m-1$ of the first $m+n-2$ spotlights are horizontal, and thus is
\begin{equation*}
\binom{m + n - 2}{m-1}.
\end{equation*}

\begin{figure}[htbp]
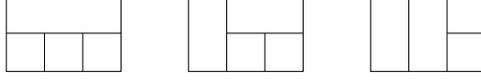

\epsfig{file=2x3-3.eps,scale=.4} \hspace{.25in} \epsfig{file=2x3-5.eps,scale=.4} \hspace{.25in} \epsfig{file=2x3-6.eps,scale=.4}\caption{The three maximal spotlight tilings of a $2 \times 3$ rectangle.  These are the spotlight tilings of Figure~\ref{fig:2x3} which contain HV-spotlights.}\label{fig:2x3-max}
\end{figure}

\begin{defn}
Let $t^r_{m,n}$ be the number of spotlight tilings of $R_{m,n}$ that use $r$ spotlights.  That is, $t^r_{m,n} = |\{t \in \mathcal{T}_{m,n} \mid |t| = r\}|$.  Set $t^r_{m,0} = t^r_{0,n} = \delta_{0r}$, where $\delta_{0r}$ is the Kronecker delta.
\end{defn}

\begin{thm}\label{thm:counting rect}
For all integers $r < m + n - 1$,
\begin{equation*}
t^r_{m,n} = \binom{r-1}{m-1} + \binom{r-1}{n-1}.
\end{equation*}
\end{thm}

Note that if $r < \max\{m,n\}$, then at least one of the binomial coefficients in the statement of the theorem is $0$, by the convention that $\binom{j}{i} = 0$ if $i > j$.

\begin{proof}

As in the proof of Theorem~\ref{thm:total rect}, the values $t^r_{m,n}$ satisfy a recurrence relation.  That is, for all $m,n,r > 0$ such that $mn > 1$,
\begin{equation*}
t^r_{m,n} = t^{r-1}_{m-1,n} + t^{r-1}_{m,n-1}.
\end{equation*}

The base case $t^1_{1,1} = 1$ is easy to calculate, and the result follows by induction.
\end{proof}

Therefore, Theorems~\ref{thm:max rect} and~\ref{thm:counting rect} and Proposition~\ref{prop:bounds} can be combined in the following equation:
\begin{equation*}
t^r_{m,n} = 
\begin{cases}
\rule[-3mm]{0mm}{7mm}\binom{r-1}{m-1} + \binom{r-1}{n-1} & r < m+n-1;\\
\binom{m+n-2}{m-1} & r = m+n-1.
\end{cases}
\end{equation*}

Observe that $t^{m+n-1}_{m,n}$ is exactly half of $\binom{m+n-1-1}{m-1} + \binom{m+n-1-1}{n-1}$, which would have been the value if Theorem~\ref{thm:counting rect} had applied.  This differences arises from the HV-spotlight present in any maximal spotlight tiling.  If the orientation of such a spotlight could be distinguished, then there would be twice as many maximal spotlight tilings of the rectangle.  As suggested earlier, the convention in this paper that an HV-spotlight lose its orientation supports the idea that these dynamic spotlight tilings should be considered as coverings of a region, and so are only distinguished if they actually look different.  However, analogously concise enumeration results will arise if this convention is dropped.

In fact, if $(m,n) \neq (1,1)$, then $t^{m+n-2}_{m,n} = t^{m+n-1}_{m,n}$, and the values $t^r_{m,n}$ are strictly increasing on the interval $r \in [\min\{m,n\}, m+n-2]$.  More specifically, for $r \in [\min\{m,n\}+1, m+n-2]$,
\begin{eqnarray*}
t^r_{m,n} - t^{r-1}_{m,n} &=& \binom{r-1}{m-1} + \binom{r-1}{n-1} - \binom{r-2}{m-1} - \binom{r-2}{n-1}\\
&=& \binom{r-2}{m-2} + \binom{r-2}{n-2} = t^{r-1}_{m-1,n-1}.
\end{eqnarray*}

\noindent Moreover, it is straightforward to check that
\begin{equation*}
\sum_{r \ge 1} t^r_{m,n} = \binom{m+n}{m} - \binom{m+n-2}{m-1},
\end{equation*}
\noindent confirming Theorem~\ref{thm:total rect}.

Given Theorems~\ref{thm:max rect} and~\ref{thm:counting rect}, it is straightforward now to compute the average number of spotlights used in a spotlight tiling of an $m \times n$ rectangle.

\begin{cor}
The average number of spotlights used in a spotlight tiling of $R_{m,n}$, that is, the average size of an element of $\mathcal{T}_{m,n}$, is
\begin{equation}\label{eqn:average}
\frac{mn(m+n-1)}{(m+n)(m+n-1)-mn}\left(1 + \frac{n-1}{m+1} + \frac{m-1}{n+1}\right).
\end{equation}
\end{cor}

\begin{proof}
This average is computed by evaluating
\begin{eqnarray*}
\frac{\sum\limits_{r=1}^{m+n-1}r\cdot t^r_{m,n}}{\binom{m+n}{m} - \binom{m+n-2}{m-1}} &=& \frac{(m+n-1)\binom{m+n-2}{m-1} + \sum\limits_{r=1}^{m+n-2}\left[r\binom{r-1}{m-1} + r\binom{r-1}{n-1}\right]}{\binom{m+n}{m} - \binom{m+n-2}{m-1}}\\
&=& \frac{(m+n-1)\binom{m+n-2}{m-1} + m\binom{m+n-1}{m+1} + n\binom{m+n-1}{n+1}}{\binom{m+n}{m} - \binom{m+n-2}{m-1}}\\
&=& \frac{mn(m+n-1)}{(m+n)(m+n-1)-mn}\left(1 + \frac{n-1}{m+1} + \frac{m-1}{n+1}\right).
\end{eqnarray*}
\end{proof}

The growth of the expression in \eqref{eqn:average} can be seen in Table~\ref{table:averages}, which displays the expected number of spotlights in a random spotlight tiling of $R_{m,n}$ for small values of $m$ and $n$.  Additionally, the average number of spotlights used in a spotlight tiling of the square $R_{n,n}$ approaches $2n-7/3$ as $n$ increases, as reflected in the table.

\begin{table}[htbp]
\centering
\begin{tabular}{c|ccccccc}
\rule[-2mm]{0mm}{6mm} & $n=1$ & $2$ & $3$ & $4$ & $5$ & $6$ & $7$\\
\hline
\rule[0mm]{0mm}{4mm}$m=1$ &$1$ & $1.5$ & $2$ & $2.5$ & $3$ & $3.5$ & $4$\\
$2$ & $1.5$ & $2.5$ & $3.286$ & $4$ & $4.688$ & $5.364$ & $6.034$\\
$3$ & $2$ & $3.286$ & $4.286$ & $5.16$ & $5.976$ & $6.762$ & $7.533$\\
$4$ & $2.5$ & $4$ & $5.16$ & $6.16$ & $7.077$ & $7.948$ & $8.793$\\
$5$ & $3$ & $4.688$ & $5.976$ & $7.077$ & $8.077$ & $9.018$ & $9.923$\\
$6$ & $3.5$ & $5.364$ & $6.762$ & $7.948$ & $9.018$ & $10.018$ & $10.974$\\
$7$ & $4$ & $6.034$ & $7.533$ & $8.793$ & $9.923$ & $10.974$ & $11.934$\\
\end{tabular}
\smallskip
\caption{The average number of spotlights used in a spotlight tiling of $R_{m,n}$ for $m,n \in [1,7]$, rounded to three decimal places.}\label{table:averages}
\end{table}

In a maximal spotlight tiling of $R_{m,n}$, there are $m-1$ horizontal spotlights, $n-1$ vertical spotlights, and $1$ HV-spotlight.  Moreover, a spotlight tiling $t \in \mathcal{T}_{m,n}$ contains an HV-spotlight if and only if $t$ is maximal.  The breakdown of spotlight directions is immediate for maximal spotlight tilings, but the question is more subtle for non-maximal spotlight tilings.

\begin{defn}
For a spotlight tiling $t$ with no HV-spotlights, let $h(t)$ be the number of horizontal spotlights in $t$, and let $v(t)$ be the number of vertical spotlights in $t$.
\end{defn}

\begin{defn}
Define the generating function
\begin{equation*}
G_{m,n}(H,V) = \sum_{\genfrac{}{}{0pt}{}{\text{non-maximal}}{t \in \mathcal{T}_{m,n}}} H^{h(t)} V^{v(t)}.
\end{equation*}
\end{defn}

Notice that $G_{1,1}(H,V) = 0$, because the only spotlight tiling of a $1 \times 1$ rectangle is maximal, yielding an empty sum.

\begin{thm}\label{thm:rect HV count}
For all $m, n \ge 1$, where $(m,n) \neq (1,1)$,
\begin{equation*}
G_{m,n}(H,V) = H^m \sum_{r=0}^{n-2} \binom{r+m-1}{m-1} V^r + V^n \sum_{r=0}^{m-2} \binom{r+n-1}{n-1} H^r.
\end{equation*}
\end{thm}

\begin{proof}

Consider a non-maximal spotlight tiling of $R_{m,n}$ using $r$ spotlights.  In the successive iterations of the spotlight tiling procedure, the last untiled sub-rectangle will be covered either by a horizontal or by a vertical spotlight.  Thus, after placing the first $r-1$ spotlights, what remains must be a rectangle of dimensions $1 \times (m + n - r)$ or $(m + n - r) \times 1$.  In the former case, the final spotlight is horizontal, and in the latter case the final spotlight is vertical.

In the case of a final horizontal spotlight, there are $m-1$ of the first $r-1$ spotlights which are horizontal, and the remaining $r - m$ are vertical.  The recursive nature of spotlight tiling means that these horizontal and vertical spotlights can occur in any order.  Thus there are $\binom{r-1}{m-1}$ ways for the last spotlight to be horizontal in a non-maximal element of $\mathcal{T}_{m,n}$ with $r$ spotlights.  Similarly, there are $\binom{r-1}{n-1}$ ways for the last spotlight to be vertical in a non-maximal element of $\mathcal{T}_{m,n}$ with $r$ spotlights.

Therefore,
\begin{eqnarray*}
G_{m,n}(H,V) &=& \sum_{\genfrac{}{}{0pt}{}{\text{non-maximal}}{t \in \mathcal{T}_{m,n}}} H^{h(t)} V^{v(t)}\\
&=& \sum_{r = \min\{m,n\}}^{m+n-2} \binom{r-1}{m-1}H^{m-1}V^{r-m}\cdot H\\
& & \hspace{.4in} + \sum_{r = \min\{m,n\}}^{m+n-2} \binom{r-1}{n-1}V^{n-1}H^{r-n}\cdot V\\
&=& \sum_{r = m}^{m+n-2} \binom{r-1}{m-1}H^mV^{r-m} + \sum_{r = n}^{m+n-2} \binom{r-1}{n-1}V^nH^{r-n}\\
&=& H^m \sum_{r=0}^{n-2} \binom{r+m-1}{m-1} V^r + V^n \sum_{r=0}^{m-2} \binom{r+n-1}{n-1} H^r.
\end{eqnarray*}
\end{proof}

One consequence of Theorem~\ref{thm:rect HV count} is that in any non-maximal spotlight tiling of $R_{m,n}$, there are either exactly $m$ horizontal spotlights or exactly $n$ vertical spotlights.  In the former case, there can be between $0$ and $n-2$ vertical spotlights, and in the latter case there can be between $0$ and $m-2$ horizontal spotlights.

Substituting $x$ for both $H$ and $V$ in $G_{m,n}(H,V)$ gives the generating function for the numbers $t^r_{m,n}$ when $r < m + n - 1$, and in fact the coefficient $[x^r]G_{m,n}(x,x)$ is equal to $\binom{r-1}{m-1} + \binom{r-1}{n-1}$, confirming Theorem~\ref{thm:counting rect}.

\section{Spotlight tilings of rectangles with missing corners}\label{sec:other regions}

The recursive nature of spotlight tilings means that enumerating the spotlight tilings of certain families of regions can be done without difficulty.  For the most part, the regions considered in this section are variations on rectangles, in particular rectangles missing squares at the corners.  Because the northwest corner is specified in spotlight tilings, the enumeration of the spotlight tilings of these regions depends on which corner was removed.

It should be noted that it is possible to obtain formulae for the number of spotlight tilings of other regions as well, due to the iterative definition of this method.  For example, the number of spotlight tilings of a rectangle with a single square removed from somewhere in the interior is not difficult to obtain, particularly if this square is parameterized by its position relative to the southeast corner of the rectangle, which does not change when spotlights are placed.

\begin{defn}
Fix integers $m, n \ge 2$.  Let $R_{m,n}^{\sf{NW}}$ (respectively, $R_{m,n}^{\sf{NE}}$, $R_{m,n}^{\sf{SW}}$, and $R_{m,n}^{\sf{SE}}$) be an $m \times n$ rectangle whose northwest (respectively, northeast, southwest, and southeast) corner has been removed.  The set $\mathcal{T}_{m,n}^*$ consists of all spotlight tilings of the region $R^*_{m,n}$, and $T_{m,n}^* = |\mathcal{T}_{m,n}^*|$.
\end{defn}

The most difficult of these spotlight tilings to enumerate, and the one with the least elegant answer, is for the region $R_{m,n}^{\sf{NW}}$.  That this case differs from the others is no surprise, since there are two northwest corners in the new region, and thus spotlights can start from two different squares.

\begin{prop}\label{prop:northwest corner}
For all $m, n \ge 2$,
\begin{eqnarray*}
T_{m,n}^{\sf{NW}} &=& T_{m-1,n-1} + T_{1,n-1}T_{m-2,n} + T_{m-1,1}T_{m,n-2}\\
&=& T_{m-1,n-1} + (n-1)T_{m-2,n} + (m-1)T_{m,n-2}\\
&=& \binom{m+n-2}{m-1} \left[1 + (m-1)(n-1)\left(\frac{1}{m} + \frac{1}{n} - \frac{1}{m+n-2}\right)\right]
\end{eqnarray*}
\end{prop}

Just as Proposition~\ref{prop:northwest corner} computes $T_{m,n}^{\sf{NW}}$, the spotlight tilings of $R_{m,n}^{\sf{NE}}$, $R_{m,n}^{\sf{SW}}$, and $R_{m,n}^{\sf{SE}}$ can also be enumerated.  In fact, these enumerations are significantly more elegant, due to the fact that the missing corner does not affect where spotlights may begin.  The proofs of these results are inductive, and use the recursion inherent to spotlight tilings.

\begin{prop}\label{prop:northeast/southwest corner}
For all $m, n \ge 2$, the number of spotlight tilings of an $m \times n$ rectangle missing either its northeast or its southwest corner is
\begin{eqnarray*}\label{eqn:T_{m,n}^{NE}}
T_{m,n}^{\sf{NE}} = T_{m,n}^{\sf{SW}} &=& T_{m,n} - 1\\
&=& \binom{m+n}{m} - \binom{m+n-2}{m-1} - 1.
\end{eqnarray*}
\end{prop}

\begin{prop}\label{prop:southeast corner}
For all $m, n \ge 2$, the number of spotlight tilings of an $m \times n$ rectangle missing its southeast corner is
\begin{eqnarray*}
T_{m,n}^{\sf{SE}} &=& T_{m,n} - \binom{m+n-2}{m-1}\\
&=& \binom{m+n}{m} - 2 \binom{m+n-2}{m-1}.
\end{eqnarray*}
\end{prop}

\begin{proof}
The number of spotlight tilings of $R_{m,n}^{\sf{SE}}$ is the number of spotlight tilings of $R_{m,n}$, minus the number of maximal spotlight tilings of $R_{m,n}$.
\end{proof}

The numbers described in Proposition~\ref{prop:southeast corner} are sequence A051601 in \cite{oeis}.

While the symmetry $T_{m,n}^{\sf{NE}} = T_{n,m}^{\sf{SW}}$ in Proposition~\ref{prop:northeast/southwest corner} is not surprising, the fact that $T_{m,n}^{\sf{NE}}$ (and $T_{m,n}^{\sf{SW}}$) is symmetric with respect to $m$ and $n$ is intriguing.  Similarly, the fact that the results of Propositions~\ref{prop:northeast/southwest corner} and~\ref{prop:southeast corner} are so similar to $T_{m,n}$ indicates that removing one of these corners does not drastically alter the spotlight tilings of the original rectangle.

In fact, Proposition~\ref{prop:northeast/southwest corner} could also be proved in another fashion, which highlights a more general trend in spotlight tilings.

\begin{defn}\label{defn:R[r]}
Suppose that $R$ is a region as in the following figure, where the only requirement of $R$ in the dashed portion is that it have no northwest corners there.
\begin{equation*}
\begin{picture}(0,0)%
\epsfig{file=NEwithR.pstex}%
\end{picture}%
\setlength{\unitlength}{1579sp}%
\begingroup\makeatletter\ifx\SetFigFont\undefined%
\gdef\SetFigFont#1#2#3#4#5{%
  \reset@font\fontsize{#1}{#2pt}%
  \fontfamily{#3}\fontseries{#4}\fontshape{#5}%
  \selectfont}%
\fi\endgroup%
\begin{picture}(4554,2571)(3364,-3970)
\put(7426,-2761){\makebox(0,0)[lb]{\smash{{\SetFigFont{10}{12.0}{\familydefault}{\mddefault}{\updefault}$r$}}}}
\put(5476,-1786){\makebox(0,0)[lb]{\smash{{\SetFigFont{10}{12.0}{\familydefault}{\mddefault}{\updefault}$n$}}}}
\put(3526,-2761){\makebox(0,0)[lb]{\smash{{\SetFigFont{10}{12.0}{\familydefault}{\mddefault}{\updefault}$r$}}}}
\end{picture}%

\end{equation*}
Let $R[r]$ be the region obtained from $R$ be removing the top $r$ squares in the rightmost column specified in $R$.  That is, $R[r]$ is the region displayed below.
\begin{equation*}
\begin{picture}(0,0)%
\epsfig{file=NEwithoutR.pstex}%
\end{picture}%
\setlength{\unitlength}{1579sp}%
\begingroup\makeatletter\ifx\SetFigFont\undefined%
\gdef\SetFigFont#1#2#3#4#5{%
  \reset@font\fontsize{#1}{#2pt}%
  \fontfamily{#3}\fontseries{#4}\fontshape{#5}%
  \selectfont}%
\fi\endgroup%
\begin{picture}(4179,2571)(3364,-3970)
\put(5026,-1786){\makebox(0,0)[lb]{\smash{{\SetFigFont{10}{12.0}{\familydefault}{\mddefault}{\updefault}$n-1$}}}}
\put(3526,-2761){\makebox(0,0)[lb]{\smash{{\SetFigFont{10}{12.0}{\familydefault}{\mddefault}{\updefault}$r$}}}}
\put(7051,-2761){\makebox(0,0)[lb]{\smash{{\SetFigFont{10}{12.0}{\familydefault}{\mddefault}{\updefault}$r$}}}}
\end{picture}%

\end{equation*}
The column of $r$ squares which gets removed from $R$ to form $R[r]$ is the \emph{difference column}.
\end{defn}

By this definition, $R_{m,n}^{\sf{NE}} = R_{m,n}[1]$.

\begin{prop}\label{prop:corner column}
Let $R$ and $R[r]$ be regions defined as in Definition~\ref{defn:R[r]}, keeping the meaning of $r$ and $n$.  Then
\begin{equation*}
\#\{\text{spotlight tilings of } R[r]\} = \#\{\text{spotlight tilings of } R\} - \sum_{k=0}^{r-1} \binom{n-1}{k}.
\end{equation*}
\end{prop}

\begin{proof}
Consider the ways that the difference column might be tiled by spotlights in $R$.  It can consist of the ends of $r$ horizontal spotlights, or the ends of $k$ horizontal spotlights atop a vertical spotlight, where $0 \le k \le r-1$.
If a vertical spotlight is involved, then this spotlight would continue down below the difference column into $R[r] \subset R$.  Additionally, if a vertical spotlight is used to cover the difference column, then there must be $n-1$ other vertical spotlight tiles positioned to the left of the difference column in $R$.  The placement of these $n-1$ vertical spotlight tiles and the $k$ horizontal spotlight tiles can be done in any order.

A given spotlight tiling of $R[r]$ can be extended to a spotlight tiling of $R$ by filling the difference column with horizontal spotlights (if the spotlight tiling of $R[r]$ included a horizontal terminating at the difference column in some row, then glue an extra square to the end of this spotlight tile).  This will yield all spotlight tilings of $R$ except those which cover some portion of the difference column with a vertical spotlight tile.  This concludes the proof.
\end{proof}

Notice that Proposition~\ref{prop:corner column} agrees with Proposition~\ref{prop:northeast/southwest corner} by setting $r = 1$.

Also notice that the symmetry of spotlight tilings indicates that Proposition~\ref{prop:corner column} would also be true if the figures in Definition~\ref{defn:R[r]} were reflected across the northwest-southeast diagonal.

One specific corollary to Proposition~\ref{prop:corner column} is presented below, although this could also have been shown in a straightforward proof using the recursion inherent to spotlight tilings.

\begin{defn}
Fix integers $m, n \ge 3$.  Let $R_{m,n}^{\sf{NE,SE}}$ be the region obtained from $R_{m,n}$ by removing the northeast and southeast corners.  Likewise, $R_{m,n}^{\sf{NE,SW,SE}}$ is an $m \times n$ rectangle whose northeast, southwest, and southeast corners have been removed.  Other regions are defined analogously, and $\mathcal{T}_{m,n}^*$ and $T_{m,n}^*$ have their customary definitions.
\end{defn}

\begin{cor}
For all $m, n \ge 3$
\begin{eqnarray*}
T_{m,n}^{\sf{NE,SW}} &=& T_{m,n} - 2\\
&=& \binom{m+n}{m} - \binom{m+n-2}{m-1} - 2;\\
\\
T_{m,n}^{\sf{NE,SE}} = T_{m,n}^{\sf{SW,SE}} &=& T_{m,n}^{\sf{SE}} - 1\\
&=& \binom{m+n}{m} - 2\binom{m+n-2}{m-1} - 1;\\
\\
T_{m,n}^{\sf{NE,SW,SE}} &=& T_{m,n}^{\sf{SE}} - 2\\
&=& \binom{m+n}{m} - 2\binom{m+n-2}{m-1} - 2.
\end{eqnarray*}
\end{cor}

There are several regions $R_{m,n}^*$ whose spotlight tilings have not yet been enumerated.  In these, the northwest corner has been removed, along with at at least one other corner.  Six of these seven cases are treated in Corollary~\ref{cor:missing corners}, and the remaining case (when all four corners have been removed) appears independently below.  The results of Corollary~\ref{cor:missing corners} are not written in closed form, although it would not be hard to do so.

\begin{cor}\label{cor:missing corners}
For $m, n \ge 3$,
\begin{eqnarray*}
T_{m,n}^{\sf{NW,SE}} &=& T_{m-1,n-1}^{\sf{SE}} + (n-1)T_{m-2,n}^{\sf{SE}} + (m-1)T_{m,n-2}^{\sf{SE}};\\
\\
T_{m,n}^{\sf{NW,NE}} &= &T_{n,m}^{\sf{NW,SW}}\\
 &=& T_{m-1,n-1} + (n-2)T_{m-2,n} +(m-1)T_{m,n-2} - m + 1;\\
\\
T_{m,n}^{\sf{NW,NE,SE}} &=& T_{n,m}^{\sf{NW,SW,SE}}\\
 &=& T_{m-1,n-1}^{\sf{SE}} + (n-2)T_{m-2,n}^{\sf{SE}} + (m-1)T_{m,n-2}^{\sf{SE}} - m + 1;\\
\\
T_{m,n}^{\sf{NW,NE,SW}} &=& T_{m-1,n-1} + (n-2)T_{m-2,n} + (m-2)T_{m,n-2} - m - n + 4.
\end{eqnarray*}
\end{cor}

\begin{defn}
For $m, n \ge 3$, let $R_{m,n}^{\circ}$ be the region obtained from $R_{m,n}$ by removing the northwest, northeast, southwest, and southeast corner squares.  Let $\mathcal{T}_{m,n}^{\circ}$ be the set of spotlight tilings of $R_{m,n}^{\circ}$, and $T_{m,n}^{\circ} = |\mathcal{T}_{m,n}^{\circ}|$.
\end{defn}

The following formula for $T_{m,n}^{\circ}$ is not difficult to compute, using the inductive definition of spotlight tilings.

\begin{cor}\label{cor:missing all corners}
For all $m,n \ge 3$,
\begin{equation*}
T_{m,n}^{\circ} = T_{m-1,n-1}^{\sf SE} + (n-2)T_{m-2,n}^{\sf SE} + (m-2)T_{m,n-2}^{\sf SE} - m - n + 4.
\end{equation*}
\end{cor}

The similarities between the results in Corollaries~\ref{cor:missing corners} and~\ref{cor:missing all corners} are striking, and suggest that the iterative nature of spotlight tiling respects certain substructures of a region.

\section{Spotlight tilings of frame-like regions}\label{sec:frames}

This section explores the spotlight tilings of a family of regions that are formed by making a large hole in the center of a rectangle.  To give a flavor for these results, this discussion studies only those cases where the remaining region has width $1$, although it would not be difficult to generalize to wider frames.

\begin{defn}
Fix $m, n \ge 3$.  Let $F_{m,n}$ be the region formed by removing a centered $(m-2) \times (n-2)$ rectangle from the rectangle $R_{m,n}$.  Let $f_{m,n}$ be the number of spotlight tilings of $F_{m,n}$.
\end{defn}

In other words, the region $F_{m,n}$ looks like an $m \times n$ picture frame of width $1$.  To understand $f_{m,n}$, it is helpful first to enumerate the spotlight tilings of some related regions.

\begin{defn}
Fix $m, n \ge 1$.  Let $C_{m,n}^{\sf{NW}}$ be the region of $m + n -1$ squares formed by overlapping the north-most square of a column of length $m$ and the west-most square of a row of length $n$.  Let $c_{m,n}^{\sf{NW}}$ be the number of spotlight tilings of $C_{m,n}^{\sf{NW}}$. The regions $C_{m,n}^{\sf{NE}}$, $C_{m,n}^{\sf{SW}}$, and $C_{m,n}^{\sf{SE}}$ and their enumerations are defined analogously.
\end{defn}

\begin{prop}\label{prop:corners}
For $m, n \ge 1$,
\begin{eqnarray*}
c_{m,n}^{\sf{NW}} &=& m + n -2\\
c_{m,n}^{\sf{NE}} = c_{n,m}^{\sf{SW}} &=& n(m-1) + 1\\
c_{m,n}^{\sf{SE}} &=& 2(m-1)(n-1) + 1
\end{eqnarray*}
\end{prop}

\begin{proof}
Each of these quantities can be computed by careful counting, together with the fact that $T_{1,p} = T_{p,1} = p$.
\end{proof}

\begin{thm}\label{thm:frames}
For $m, n \ge 3$,
\begin{equation*}
f_{m,n} = 2(m-2)(n-2)(m+n-2) + (m-2)(m+1) + (n-2)(n+1).
\end{equation*}
\end{thm}

\begin{proof}
Initially, there is only one northwest corner in the region $F_{m,n}$.  This can be covered with a horizontal spotlight of length $n$ or a vertical spotlight of length $m$.  Either way, the remaining region has two northwest corners, and careful applications of Proposition~\ref{prop:corners} and the inclusion-exclusion property give the answer.
\end{proof}

The values of $f_{m,n}$ for small $m$ and $n$ are displayed in Table~\ref{table:f_{m,n}}.  These values are sequence A132370 of \cite{oeis}.

\begin{table}[htbp]
\centering
\begin{tabular}{c|ccccc}
\rule[-2mm]{0mm}{6mm}$f_{m,n}$ & $n=3$ & 4 & 5 & 6 & 7\\
\hline
\rule[0mm]{0mm}{4mm}$m=3$ & 16 & 34 & 58 & 88 & 124\\
4 & 34 & 68 & 112 & 166 & 230\\
5 & 58 & 112 & 180 & 262 & 358\\
6 & 88 & 166 & 262 & 376 & 508\\
7 & 124 & 230 & 358 & 508 & 680
\end{tabular}
\smallskip
\caption{The number of spotlight tilings of $F_{m,n}$, for $m,n \in [1,7]$.}\label{table:f_{m,n}}
\end{table}

\section{Further directions}\label{sec:further}

The preceding sections have examined the spotlight tilings of several families of regions.  In each case, the enumeration of these spotlight tilings had a concise and often illuminating form.  For the rectangle, more refined analysis was also performed, and yielded results whose simplicity and elegance may not have been anticipated.

The obvious analogue of spotlight tiling in higher dimensions may also yield fruitful results.  Additionally, the questions particular to spotlight tiling (such as the distribution of the number of spotlights in a given spotlight tiling) may give rise to new aspects of this and other tilings methods which warrant further study.

This work can be extended by studying variations on the spotlight tilings described here.  For example, in this article, every spotlight has started in a northwest corner.  If this restriction were removed, and spotlights were allowed to start in any square and continue in any direction until reaching a barrier, then the resulting model would be an expansion of this type of dynamic tiling.

Another generalization would be to allow tiles to expand as much as possible in two directions, instead of only horizontally or only vertically.  Such a tile would create an $a \times b$ rectangle, instead of $a \times 1$ or $1 \times b$.  Continuing the imagery of this article, these new tiles could be called \emph{floodlights}, and dynamic floodlight tiling might have interesting enumerative results as well.  It should be noted that the region $R_{m,n}$ has exactly $1$ floodlight tiling, and, consequently, more complicated regions need to be studied in order to gain an understanding of the model.

\end{document}